\documentclass[11pt, reqno]{amsart}
\pdfoutput=1
\usepackage{amscd, amsfonts, amssymb, pinlabel}
\usepackage{graphicx}
\usepackage{xfrac}
\usepackage[breaklinks=true, backref=page]{hyperref}

\newtheorem{thm}{Theorem}

\theoremstyle{definition}

\theoremstyle{remark}

\newtheorem{rmk}[thm]{Remark}

\newtheorem{ex}[thm]{Example}

\numberwithin{thm}{section}

\numberwithin{equation}{section}

\def\R{\mathbb{R}}
\def\Z{\mathbb{Z}}

\newcommand{\C}{{\mathbb{C}}}

\newcommand{\pa}{\partial}

\newcommand{\tr}{\operatorname{tr}}

\begin{document}

\title{Notes on topological strings and knot contact homology}

\author{Tobias Ekholm}
\address{Uppsala University, Box 480, 751 06 Uppsala, Sweden\newline
\indent Insitute Mittag-Leffler, Aurav 17, 182 60 Djursholm, Sweden}
\email{tobias.ekholm@math.uu.se}
\thanks{The author was partially funded by the Knut and Alice Wallenberg Foundation as a Wallenberg scholar and by the Swedish Research Council, 2012-2365.}
\begin{abstract}
We give an introduction to the physics and mathematics involved in the recently observed relation between topological string theory and knot contact homology and then discuss this relation. The note is based on two lectures given at the G\"okova Geometry and Topology Conference, 2013 and reports on joint work by Aganagic, Ng, Vafa, and the author \cite{AENV}.   
\end{abstract}

\maketitle


\section{Introduction}\label{Sec:intro}
These notes introduce,  discuss, and provide background for the recently observed relations between topological string theory and knot contact homology found in \cite{AENV}. The starting point for these relations was the mirror symmetry motivated construction of Aganagic and Vafa \cite{AV} of a polynomial associated to any knot. In concrete calculations, this new polynomial was observed to agree with Ng's augmentation polynomial of knot contact homology, see \cite{Ngframed, Ngtransverse}, and the two polynomials are conjectured to agree for any knot. Here we will explain how to relate the algebraic varieties determined by the respective polynomials and we will show that certain branches of them agree. In fact, our argument proves that the polynomials agree if the corresponding curves are known to be irreducible (but no general condition on knots that guarantees irreducibility is known). In \cite{AENV}, an analogue of this relation for many component links is presented. That subject will be discussed only very briefly in these notes, see Section \ref{s:links}.  

An outline of the paper is as follows. In Section \ref{S:CSandtopstring} we discuss necessary background on Chern-Simons theory and topological strings from a physics perspective. In Section \ref{S:largeN} we recall the relation between topological strings in the cotangent bundle of the three-sphere and in the resolved conifold and describe the mirror construction from \cite{AV}. In Section \ref{S:kch} we introduce knot contact homology, and in Section \ref{S:aug} we define the augmentation variety and demonstrate how it is related to the polynomial in the mirror construction.
     
\subsection*{Acknowledgments} The author thanks Mina Aganagic, Lenny Ng, and Cumrun Vafa for very useful comments on the text.

\section{Chern-Simons theory and topological strings}\label{S:CSandtopstring}
We give a physics description of Chern-Simons theory, topological strings, and their relation. 
Much of the material here is originally due to Witten. Except for original references e.g.~ \cite{Witten1986253, Witten_top_sigma, Witten:1988hf, Witten:1991zz, Witten:1992fb}, the presentation is also based on \cite{Marino}. 

\subsection{Chern-Simons action}
Let $A$ be an $U(N)$-connection on a closed 3-manifold $M$, i.e.~a 1-form on $M$ with values in complex $N\times N$-matrices $A$ with $A=-A^{\ast}$. Then the Chern-Simons action of $A$ is the integral
\[
S(A)=\int_{M} \tr\left(A\wedge dA + \tfrac{2}{3}A\wedge A\wedge A\right).
\]

Consider a variation $A\mapsto A+\delta A$ of $A$. The corresponding variation of $S$ is 
\[
\delta S = 2\int_{M} \tr\left(\delta A\wedge F_{A}\right),
\]
where $F_{A}=dA+A\wedge A$ is the curvature of $A$. To see this note that
\[
\int_{M}\delta A\wedge dA + A\wedge d(\delta A) = 2\int_{M}\delta A\wedge dA,
\]
by Stokes theorem. Furthermore, writing $A=(A^{i}_{j})$ and using the summation convention, 
$\tr(A\wedge A\wedge A)=A^{i}_{j}\wedge A^{j}_{k}\wedge A^{k}_{i}$:
\begin{align}\label{eq:varqube}
&\delta A^{i}_{j}\wedge A^{j}_{k}\wedge A^{k}_{i}+
A^{i}_{j}\wedge \delta A^{j}_{k}\wedge A^{k}_{i}+
A^{i}_{j}\wedge A^{j}_{k}\wedge \delta A^{k}_{i}=\\\notag
&
\delta A^{i}_{j}\wedge A^{j}_{k}\wedge A^{k}_{i}
-\delta A^{j}_{k}\wedge A^{i}_{j}\wedge A^{k}_{i}+
\delta A^{k}_{i}\wedge A^{i}_{j}\wedge A^{j}_{k}=\\\notag
&\delta A^{i}_{j}\wedge A^{j}_{k}\wedge A^{k}_{i}
+\delta A^{j}_{k}\wedge A^{k}_{i}\wedge A^{i}_{j}+
\delta A^{k}_{i}\wedge A^{i}_{j}\wedge A^{j}_{k}=
3\;\delta A^{i}_{j}\wedge A^{j}_{k}\wedge A^{k}_{i}.
\end{align}
In conclusion, connections that are stationary for $S$ are flat.

We next consider the effect of a gauge transformation, $g\colon M\to U(N)$. The covariant derivative is $d+A$. Changing trivialization with $g$, we find that the covariant derivative in the new trivialization is 
\[
g^{-1}(d+A)g=d+(g^{-1}Ag+g^{-1}dg),
\]   
and thus the connection transforms as $A\mapsto A'=g^{-1}A g+g^{-1}dg$. The Chern-Simons action is not quite invariant under gauge transformations, rather
\begin{equation}\label{eq:intjump}
S(A')=S(A)+8\pi^{2} k,\quad k\in\Z.
\end{equation}
To see this we consider the the 4-manifold $M\times [0,1]$ and a trivial $\C^{N}$-bundle on it. We build the bundle $\xi$ on $M\times S^{1}=M\times [0,1]/(M\times\{0\}\sim M\times\{1\})$ by identifying the fibers near $M\times \{0\}$ with the fibers near $M\times\{1\}$ by the gauge transformation $g$. Pick a function $\phi\colon [0,1]\to [0,1]$ which equals $0$ near $0$ and equals $1$ near $1$. Then the connection 1-form
\[
\mathbf{A}=(1-\phi(t))A+\phi(t)A',
\] 
on $M\times[0,1]$ descends to a connection 1-form on $M\times S^{1}$. Furthermore, 
\[
d\tr(\mathbf{A}\wedge d\mathbf{A}+\tfrac{2}{3}\mathbf{A}\wedge\mathbf{A}\wedge\mathbf{A})=
\tr(F_{\mathbf{A}}^{2}).
\] 
To see this note that the calculation of the left hand side is similar to \eqref{eq:varqube} and that in the right hand side the term $\tr(\mathbf{A}^{\wedge 4})$ in $\tr(d\mathbf{A}+\mathbf{A}\wedge\mathbf{A})^{2}$ gives no contribution since
\[
\mathbf{A}^{i}_{j}\wedge\mathbf{A}^{j}_{k}\wedge\mathbf{A}^{k}_{l}\wedge\mathbf{A}^{l}_{i}=
-\mathbf{A}^{l}_{i}\wedge\mathbf{A}^{i}_{j}\wedge\mathbf{A}^{j}_{k}\wedge\mathbf{A}^{k}_{l}.
\]
Now the Chern-Weil formula for the second Chern class reads
\[
c_2(\xi)=\frac{1}{8\pi^{2}}\int_{M\times S^{1}}\tr(F_{\mathbf{A}}^{2})
\]
and 
\begin{align*}
&\int_{M\times S^{1}}\tr(F_{\mathbf{A}}^{2})=\int_{M\times[0,1]}\tr(F_{\mathbf{A}}^{2})
=\int_{M\times[0,1]}d\tr(\mathbf{A}\wedge d\mathbf{A}+\tfrac{2}{3}\mathbf{A}\wedge\mathbf{A}\wedge\mathbf{A})\\
&\quad=\int_{M}\tr\left(A'\wedge dA'-\tfrac23 A'\wedge A'\wedge A'\right)
-\int_{M}\tr\left(A\wedge dA-\tfrac23 A\wedge A\wedge A\right)\\
&\quad=S(A')-S(A),
\end{align*}
by Stokes theorem. Since $c_2(\xi)$ is an integer we conclude that \eqref{eq:intjump} holds.

\subsection{Chern-Simons partition function}
The \emph{Chern-Simons partition function} is the path integral
\[
Z_{\mathrm{CS}}(M)=\int\mathcal{D}A\; e^{\frac{ik}{4\pi}S(A)},\quad k\text{ integer}
\] 
over gauge orbits of connections. (Note that the integrand is gauge invariant.) Here $k$ is known as the level and is analogous to $\frac{1}{\hbar}$ in quantum mechanical path integrals. 

As for many other path integrals, it is not straightforward to give a mathematically satisfactory definition the Chern-Simons partition function. Much work has been done in that direction but we will not discuss that any further here. Rather we will study the partition function with physics methods.

First, in analogy with the stationary phase approximations for finite dimensional integrals, we expect the partition function to be expressed as a sum of perturbation expansions near each critical point of $S(A)$, i.e.~by the above remark, near each gauge orbit of a flat connection. The perturbation expansion is in powers of $\frac{1}{k}$ and is obtained via Feynman diagrams. For readers not familiar with Feynman calculus we start with a finite dimensional analogue.

Let $Q=Q_{ij}$ denote a positive definite quadratic form in $n$ variables and let $C=C_{ijk}$ be a trilinear form. Consider the integral (as above we use the summation convention)
\[
\int_{\R^{n}} e^{\;-\tfrac{1}{2}Q_{ij}x_ix_j + \hbar C_{ijk}x_ix_jx_k} \,dx.
\]
We will expand it in powers of $\hbar$. To this end we first expand the second factor in the integrand:
\begin{align*}
\int_{\R^{n}} &e^{\;-\tfrac{1}{2}Q_{ij}x_ix_j + \hbar C_{ijk}x_ix_jx_k} \,d=
\sum_{m=0}^{\infty}\;\;\int_{\R^{n}} e^{\;-\tfrac{1}{2}Q_{ij}x_ix_j}
\frac{1}{m!}(C_{ijk}x_ix_jx_k)^{m}\,dx\\
&=\sum_{m=0}^{\infty}\frac{1}{m!}\left.\left(\hbar C_{ijk}\frac{\partial}{\partial b_i}\frac{\partial}{\partial b_j}\frac{\partial}{\partial b_k}\right)^{m}\right|_{b=0}\int_{\R^{n}} e^{\;-\tfrac{1}{2}Q_{ij}x_ix_j + b_ix_i}\,dx.
\end{align*}
The last integral is Gaussian and it is straightforward to compute it explicitly via a change of variables: 
\begin{align*}
\int_{\R^{n}} & e^{\;-\tfrac{1}{2}Q_{ij}x_ix_j + \hbar C_{ijk}x_ix_jx_k} \,dx=\\
&\left(\frac{\det A}{2\pi}\right)^{-\frac{n}{2}}\sum_{m=0}^{\infty}\frac{1}{m!}\left.\left(\hbar C_{ijk}\frac{\partial}{\partial b_i}\frac{\partial}{\partial b_j}\frac{\partial}{\partial b_k}\right)^{m}\right|_{b=0}\;e^{\tfrac{1}{2}Q^{ij}b_ib_j},
\end{align*}
where $Q^{ij}$ is the inverse matrix of $Q_{ij}$.

\begin{figure}[ht]
\labellist
\small
\pinlabel $C_{ijk}$ at 235 580
\pinlabel $i$ at 142 650
\pinlabel $j$ at 142 440
\pinlabel $k$ at 345 575
\endlabellist
\centering
\includegraphics[width=.7\linewidth]{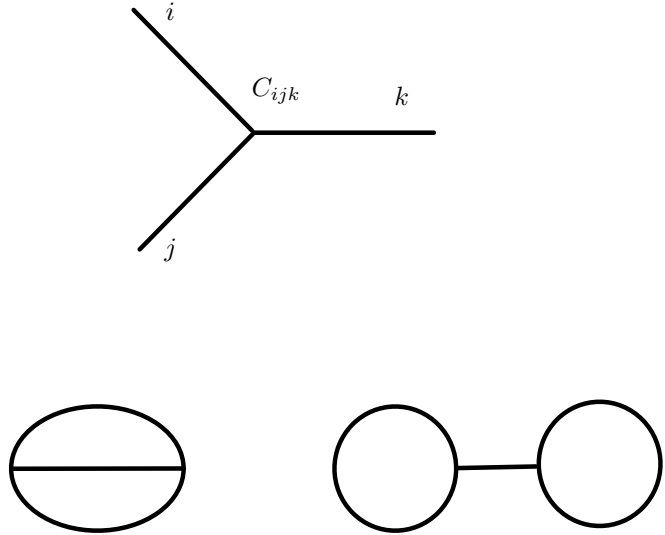}
\caption{Top: the trivalent vertex corresponding to the operator $C_{ijk}\frac{\partial}{\partial b_i}\frac{\partial}{\partial b_j}\frac{\partial}{\partial b_k}$. Bottom: the two graphs contributing to the $\hbar^{2}$-term in the expansion. }
\label{fig:feynman}
\end{figure}

We now explain how to express the right hand side as a sum over certain graphs. To this end we depict the operator $C_{ijk}\frac{\partial}{\partial b_i}\frac{\partial}{\partial b_j}\frac{\partial}{\partial b_k}$ as a trivalent vertex weighted by $C_{ijk}$, see Figure \ref{fig:feynman}. Note next that non-zero contributions to the expansion come from pairs of differential operators $\frac{\partial}{\partial b_i}$ and $\frac{\partial}{\partial b_j}$ where the first operator brings down the factor $\tfrac{1}{2}Q^{ij}b_j$ and the second removes $b_j$, here one can reverse the roles of $i$ and $j$ giving $Q^{ij}$ in total. To keep track of such contributions we join corresponding edges of the operator graphs and weight such an edge connecting $i$ to $j$ by $Q^{ij}$. With this observed we find that the contributions correspond to closed trivalent graphs with weights at trivalent vertices and edges as just explained. Furthermore, the power of $\hbar$ of such a connected graph $\Gamma$ equals its number of trivalent vertices, see Figure \ref{fig:feynman}. 
On the other hand the Euler characteristic $\chi(\Gamma)$ is, if $v$ and $e$ denote the number of vertices and edges of $\Gamma$, respectively,
\[
\chi(\Gamma)=v-e=v-\tfrac{3}{2}v=-\tfrac{1}{2}v=-\dim(H_{1}(\Gamma))+1, 
\] 
and thus the power of $\hbar$ corresponding to $\Gamma$ is $\hbar^{2(r-1)}$, where $r$ is the number of loops in $\Gamma$. Finally, we remark that the sum over all disconnected graphs equals the exponential of the sum over all connected graphs. 

\begin{figure}[ht]
\labellist
\small
\pinlabel $i$ at 207 570
\pinlabel $j$ at 300 615
\pinlabel $k$ at 300 525
\pinlabel $i$ at 60 405
\pinlabel $j$ at 90 330
\pinlabel $k$ at 190 330
\pinlabel $i$ at 445 360
\pinlabel $j$ at 345 330
\pinlabel $k$ at 560 330
\endlabellist
\centering
\includegraphics[width=.7\linewidth]{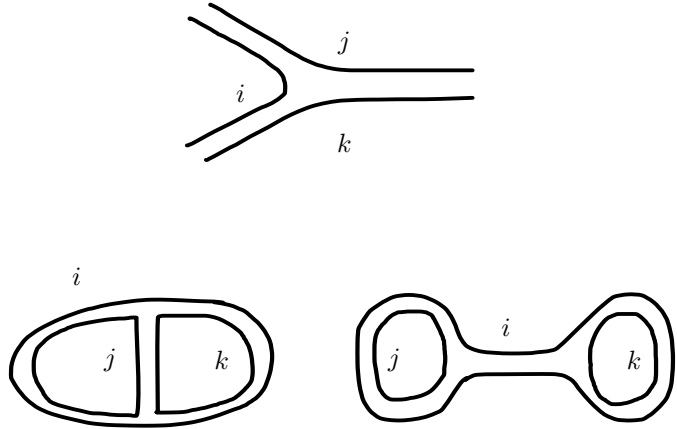}
\caption{Fattened vertex and diagrams for the $\hbar^{2}$-term. Boundary components are decorated with gauge indices.}
\label{fig:fatgraph}
\end{figure}

We now return to the infinite dimensional setting and treat the contribution $Z^{(C)}_{\mathrm{CS}}(M)$ to the partition function $Z_{\mathrm{CS}}(M)$ at a flat connection $C$ as if the path integral were a finite dimensional integral as above. Here the first term $\tr(A\wedge dA)$ in the Chern-Simons action is analogous to quadratic form $Q_{ij}x_ix_j$ above and the qubic term $\frac{2}{3}\tr(A\wedge A\wedge A)$ is analogous to the trilinear form $C_{ijk}x_ix_jx_k$.  We find that the partition function should be expressed as a sum over diagrams with trivalent vertices. In fact, the diagrams are independent of the gauge indices and to keep track of the contributions it is convenient to fatten the graphs to a thin surface, a so called a fat-graph, see Figure \ref{fig:fatgraph}.
We thus get an expansion for the free energy, i.e.~the logarithm of the partition function, of the form
\[
F_{\mathrm{CS}}^{(C)}(M)=\sum_{h,r} C'_{h,r} N^{h} k^{1-r},
\] 
where we sum over connected fat graphs with $h$ boundary components and $r$ loops, and where $F^{(C)}_{\mathrm{CS}}$ is the free energy corresponding to the local contribution $Z_{\mathrm{CS}}^{(C)}$ to the partition function from the flat connection $C$. The factor $N^{h}$ arises from the $N$ possible choices of gauge index on each boundary component. One should however rather use the effective expansion parameter or the shifted coupling constant, see e.g.~\cite{Marino} for an explanation, $\kappa=\frac{2\pi i}{k+N}$ and rewrite the above as
\[
F_{\mathrm{CS}}^{(C)}(M)=\sum_{h,r} C_{h,r} \,N^{h} \kappa^{1-r}.
\]
Viewing the fat-graph corresponding to a Feynman diagram $\Gamma$ as a Riemann surface $\Sigma$ of genus $g$ with $h$ boundary components we find that
\[
2-2g-h=\chi(\Sigma)=\chi(\Gamma)=1-r,
\] 
so the above expansion becomes
\[
F_{\mathrm{CS}}^{(C)}(M)=\sum_{g,h} C_{g,h} \,N^{h} \kappa^{2g+h-2}.
\] 
 
\subsection{A-model topological string} 
Our next goal is to relate the Chern-Simons partition function to the partition function for an A-model topological string theory on $T^{\ast}M$, which is also given by a path integral. In order to describe A-model topological string we start for simplicity with a description of closed strings in a Calabi-Yau threefold. (Note that $T^{\ast}M$ equipped with an almost complex structure $J$ compatible with the standard symplectic form satisfies the Calabi-Yau condition $c_1(T^{\ast}M)=0$.) 

Let $X$ be a 6-dimensional K\"ahler manifold with symplectic form $\omega$ and equipped with an (almost) complex structure $J$, which we assume satisfies the Calabi-Yau condition $c_1(X)=0$. The almost complex structure $J$ gives a decomposition of the complexified cotangent bundle, viewed as complex valued real linear functionals on $TX$, into $(i,J)$ complex linear and anti-linear functionals, respectively. There is a corresponding decomposition of the tangent bundle:  
\[
TX\otimes\C= T^{1,0}X\oplus T^{0,1}X,
\]  
where $T^{1,0}X$ is the holomorphic tangent space, annihilated by the complex anti-linear forms, and $T^{0,1}X$ is the anti-holomorphic tangent space annihilated by the complex linear forms.   

The fields of the A-model topological string theory are triples $(\phi,\chi,\psi)$, where 
\begin{itemize}
\item $\phi\colon\Sigma\to X$ is a map of a Riemann surface into $X$,
\item $\chi\colon\Sigma\to\phi^{\ast}(TX)$ is a vector field along $\Sigma$, and
\item $\psi\colon\Sigma\to \;\phi^{\ast}(TX^{1,0})\otimes {T^{\ast}}^{0,1}\Sigma\,\oplus\, 
\phi^{\ast}(TX^{0,1})\otimes {T^{\ast}}^{1,0}\Sigma$ is a section.   
\end{itemize} 

We write the action in local coordinates. Let $z=s+it$ be a local coordinate on the Riemann surface $\Sigma$ and write $\partial_{z}=\partial_{s}-i\partial_{t}$ and $\partial_{\bar z}=\partial_{s}+i\partial_{t}$. Let $\phi^{i}$, $i=1,\dots,6$ denote the components of the map $\phi$ in coordinates $(x^{1},\dots,x^{6})$ on $X$. Let $G_{ij}dx_i\otimes dx_j$ denote the metric $\omega(\cdot,J\cdot)$. We pick also local complex trivializations $(\xi^{1},\xi^{2},\xi^{3})$ of $T^{1,0}X$ and $(\xi^{\bar{1}},\xi^{\bar{2}},\xi^{\bar{3}})$ of $T^{0,1}X$, where we let the latter be the complex conjugates of the former. In this decomposition we get for example
\[
\psi=\psi^{I}_{\bar z}\,\xi^{I}\otimes d\bar z + \psi^{\bar{I}}_{z}\,\xi^{\bar{I}}\otimes dz.
\]

Consider the metric $G_{ij}$ given by $\omega(\cdot,J\cdot)$ and let $R_{ijkl}$ denote its curvature tensor. Note that the metric induces Hermitian forms $G_{I\bar{J}}$ and $G_{\bar{I}J}$ on $T^{(1,0)}X$ and $T^{(0,1)}X$, respectively, and likewise $R_{ijkl}$ induces a tensor with components $R_{I\bar{J}K\bar{L}}$. The Lagrangian for the A-model is now given by, with $d^{2}z=|dz\wedge d\bar{z}|$,
\begin{align*}
&L=\\
&2t\int_{\Sigma} d^{2}z
\left(\tfrac12 G_{ij}\partial_{z}\phi^{i}\partial_{\bar z}\phi^j
+iG_{I\bar{J}}\psi^{I}_{\bar{z}}D_{z}\chi^{\bar J}+
iG_{\bar{I}J}\psi^{\bar{I}}_{z}D_{\bar z}\chi^{I}
-R_{I\bar{I}J\bar{J}}\psi^{I}_{\bar z}\psi^{\bar{I}}_{z}\chi^{J}\chi^{\bar{J}}\right).
\end{align*}
Here the covariant derivative $D_\alpha$ is
\[
D_\alpha\chi^{i}=\partial_{\alpha}\chi^{i} + \partial_{\alpha}\phi^{j}\Gamma^{i}_{jk}\chi^{k},
\]
where $\Gamma_{jk}^{i}$ are the Christofel symbols of the metric $G_{ij}$. The partition function of the theory is then given by the path integral over these maps with tensor fields along them:
\[
Z(X)=\int \mathcal{D}\Phi\, e^{L},
\]
where $\Phi$ is short hand for all the fields of the theory.

A central feature of the theory is that it possesses a BRST symmetry (or a topological charge) $\mathcal{Q}$ which squares to zero, $\mathcal{Q}^{2}=0$, and which acts on the fields as follows
\begin{align*}
\{\mathcal{Q},\phi\}&=\chi,\\
\{\mathcal{Q},\chi\}&=0,\\
\{\mathcal{Q},\psi^{I}_{\bar{z}}\}&=i\partial_{\bar{z}}\phi^{I}-\chi^{J}\Gamma^{I}_{JK}\psi^{K}_{\bar{z}},\\
\{\mathcal{Q},\psi^{\bar{I}}_{z}\}&=i\partial_{z}\phi^{\bar{I}}-\chi^{\bar{J}}\Gamma^{\bar{I}}_{\bar{J}\bar{K}}\psi^{\bar{K}}_{z}.
\end{align*}

Noting that $\mathcal{Q}$ is an odd symmetry we find that the path integral localizes on the fixed points of the symmetry. To see this, we argue as in \cite{Witten:1992fb}: outside a neighborhood of the fixed points the $\mathcal{Q}$ action is free and that allows us to introduce a global Grassmann variable $\theta$ for the $\mathcal{Q}$-symmetry here. The $\mathcal{Q}$-invariance of $L$ then implies that the integral factors in this region with one of the factors independent of $\theta$ and the other the fermionic integral $\int d\theta=0$. Thus the partition function localizes at the fixed points of the action and can be calculated from data in an arbitrarily small neighborhood of it. In the case under study, we determine the fixed points from the action of $\mathcal{Q}$ described above: from the first equation we get $\chi=0$, then the last equations give 
\[
\partial_{\bar{z}}\phi^{I}=\partial_{z}\phi^{\bar I},
\]          
or equivalently
\[
d\phi+J \,d\phi\, i=0
\]
In conclusion the A-model partition function localizes on the moduli space of holomorphic curves. (The argument just given can be compared for instance with the Bott-localization formula.)

The same result can also be shown via a ``heat-kernel argument'' as follows. The Lagrangian $L$ is $\mathcal{Q}$-exact up to a topological term:
\[
L=it\{\mathcal{Q},V\} -t\int_{\Sigma}\phi^{\ast}\omega,
\]
where 
\[
V=\int_{\Sigma} d^{2}z\, G_{I\bar{J}} \left(\psi_{\bar{z}}^{I}\partial_{z}\phi^{\bar{J}}+
\psi_{z}^{\bar{I}}\partial_{\bar{z}}\phi^{I}\right).
\]
Noting that the second term in the expression for $L$ depends only on the homology class $A$ of the map $\phi$, we find that it just contributes a factor $e^{-ta}$ to the partition function where $a=\langle[\omega],A\rangle$ is given by pairing $A$ with the cohomology class of the symplectic form.

As the operator $\mathcal{Q}$ is the differential for BRST cohomology, $\mathcal{Q}$-exact quantities have zero expectation value and we find
\[
\frac{d}{dt}\int e^{it\{\mathcal{Q},V\}}=i\int\{\mathcal{Q},V\}e^{it\{\mathcal{Q},V\}}=0
\]  
and hence (except for the understood factor from the homology class) the partition function is $t$-independent. In the limit $t\to\infty$ one finds that only holomorphic maps contribute. 

The above path integral calculations were carried out for a fixed metric $h$ on the Riemann surface $\Sigma$. However, 
since we work in a Calabi-Yau threefold $X$, the formal dimension of any holomorphic curve is zero and thus for curves of genus $g$ the formal dimension of a holomorphic curve on a genus $g\ge 2$ Riemann surface with fixed conformal structure is $-(6g-6)$ and we must take variations of the metric into account in order to get a meaningful theory for all genera. To this end one notes that the stress energy tensor $T_{\alpha\beta}$, i.e.~the variation of the action with respect to the metric on $\Sigma$ is $\mathcal{Q}$-exact: since the second term in the action depends only on the homology class we have $T_{\alpha\beta}=\{\mathcal{Q}, b_{\alpha\beta}\}$ where $b_{\alpha\beta}=t\delta V/\delta g_{\alpha\beta}$. Consider now $6g-6$ variations $\delta^{(j)}h$, $j=1,\dots,6g-6$ of the metric and define 
\[
b^{(j)}=\int_{\Sigma}dA_{h}\, (\delta h^{(j)},b)_{h},
\]
where $dA_{h}$ is the area form of $h$ and $(\cdot,\cdot)_{h}$ is the induced inner product on the second symmetric power of $T^{\ast}\Sigma$. Then one shows, see \cite{Witten:1991zz}, that
\[
\Theta(\delta h^{(1)},\dots,\delta h^{(6g-6)})=\left\langle b^{(1)}\dots \;b^{(6g-6)}\right\rangle,
\]
where the right hand side denotes the expectation value in the path integral is a closed $(6g-6)$-form on the space of metrics. This form is clearly diffeomorphism invariant and furthermore vanishes if one of the $\delta^{j}h$ is induced by an infinitesimal diffeomorphism. This means that the form descends to the quotient of the space of metrics by diffeomorphism. The fact that it is closed then means that we can define an associated form $\Theta_{g}$ on the moduli space $M_g$ of genus $g$ surfaces, by pulling back via a section, that is well-defined up to an exact form. We then define the free energy in genus $g$ by
\[
F_g=\int_{M_{g}}\Theta_g. 
\]
As above one shows that this partition function concentrates on holomorphic curves as well. This time, we have taken variations of the metric into account and since we are working on a Calabi-Yau threefold all homomorphic curves are of formal dimension $0$. Thus after employing a suitable perturbation scheme we find that the integrals correspond to a curve count, exactly as in Gromov-Witten theory. Define the free energy $F_{\mathrm{GW}}(X)$ as the weighted sum of genus $g$ free energies over all topologies:
\[
F_{\mathrm{GW}}(X)=\sum_{g,A} C_{g,A} \,g_{s}^{2-2g}Q^{A},
\]
where $g_s$ is the string coupling constant, $C_{g,A}$ is the number of curve configurations of genus $g$ in the homology class $A\cdot\omega$, see Section \ref{s:nonexact} for a more precise definition in the case of disks. (To motivate the appearance of the string coupling constant: in the path integral, each emission or absorption of a closed string gives rise to a factor $g_s=e^{\Phi_0}$, where $\Phi_0$ is the energy of the dilaton field, see e.g.~\cite[p.~83]{Polchinski}.) The Gromov-Witten partition function is then $Z_{\mathrm{GW}}(X)=\exp(F_{\mathrm{GW}}(X))$.

\subsection{A-model open topological string in $T^{\ast}M$}\label{s:AmodelTM}
We next consider the corresponding open string theory in $T^{\ast}M$ with $N$ Lagrangian branes on the zero-section $M\subset T^{\ast}M$. The theory is completely analogous to what was considered above, we will however replace the closed Riemann surfaces above with Riemann surfaces $\Sigma$ with non-empty boundary. We write $h$ for the number of boundary components of the Riemann surface and think of it as a closed Riemann surface of genus $g$ with $h$ disks removed. Naturally we need to specify boundary conditions for our fields along $\pa\Sigma$. First we require that $\phi$ takes any boundary component into (any one of the $N$ copies of) the zero-section. Locally, this corresponds to Dirichlet conditions in the directions normal to the Lagrangian. In the remaining directions along the Lagrangian we impose Neumann conditions. The boundary conditions on the fields $\psi$ and $\chi$ are analogous: they are required to take values in $\phi^{\ast}(TL)$ on the boundary and to satisfy Neumann conditions in these directions. Open string theory also allows for emission and absorption of open strings that each give a factor $\sqrt{g_{s}}$ in the  path integral. As in the closed string case, the relevant path integral localizes on holomorphic curves and we get a partition function $Z_{\mathrm{GW}}(T^{\ast}M,N\cdot M)$ with free energy corresponding to the open Gromov-Witten potential:
\[
F_{\mathrm{GW}}(T^{\ast}M,N\cdot M)=\sum_{g,h} K_{g,h}\, N^{h}g_{s}^{2g+h-2},
\]    
where $K_{g,h}$ is the number of curves of genus $g$ with $h$ boundary components on the zero-section. The factor $N^{h}$ arises from the choice of $N$ branes for each boundary component.

\subsection{String field theory, topological string, and Chern-Simons}
In the previous sections we found two distinct physics ways of associating topological invariants to a closed 3-manifold in the form of partition functions. In this section we discuss Witten's proof that they are in fact the same. The proof goes via open string field theory with qubic interaction term.

Consider the open string theory in $T^{\ast}M$ with $N$ branes on the 0-section as described above. A string field is a functional $\mathcal{A}$ on the configuration space for open strings. This configuration space is the space of paths with endpoints on $M$, the Riemann surfaces wih boundary considered in Section \ref{s:AmodelTM} are then the corresponding world sheets that the strings trace out as they move. Witten \cite{Witten:1988hf} introduces two operations on the space of string functionals defined via path integrals as follows. Choose  a mid point on the interval parameterizing the open string. This subdivides any string $\xi$ into two halves $\xi_{L}$ and $\xi_{R}$ and then the two operations admit the following description: 
\begin{itemize}
\item Integration $\int\mathcal{A}$ defined by
\[
\int \mathcal{A} = \int \mathcal{D}\xi \;\delta(\xi_{L},\xi_{R})\mathcal{A}(\xi),
\]
where $\delta$ is the Dirac-delta on the space of strings, intuitively, the product over the Dirac deltas at all points along the string.
\item A $\star$-product defined analogously by
\begin{align*}
&\int \mathcal{A}_1\star\dots\star\mathcal{A}_{m} =\\
&\int\mathcal{D}\xi^{1}\dots\mathcal{D}\xi^{m}\; \delta(\xi_{R}^{1},\xi_{L}^{2})\dots\delta(\xi^{m-1}_R,\xi^{m}_{L})\delta(\xi_{R}^{m},\xi^{1}_L)\;\mathcal{A}_{1}(\xi^{1})\dots\mathcal{A}_{m}(\xi^{m}).
\end{align*}
\end{itemize}

The string theory has an underlying parameterization symmetry. Let $Q$ denote the corresponding BRST-operator and define the (Chern-Simons like) string field action as
\[
S=\frac{1}{2}\int \mathcal{A}\star Q\mathcal{A} + \tfrac23 \mathcal{A}\star\mathcal{A}\star\mathcal{A}.
\]
This action holds all the information of the space time dynamics of the open topological string. In a sense it corresponds to a Hamiltoninan description of the open string. 

The next step is to analyze this in the case of topological strings, we refer to \cite{Witten:1992fb} for more details. One first notes that the string fields can be expressed in terms of the expansion of the string states around its zero-mode. Here the zero-modes are $(\phi,\chi,\psi)$ where $\phi^{i}=x^{i}$ is a constant map into $M$, $\psi^{i}=\partial_{x^{i}}$ the corresponding tangent vector, and $\chi^{is}=-dx^{i}$.  We thus find, that the zero-modes have the form  
\[
\mathcal{A}_{0}=c(x)+\chi^{i}A_{i}(x)+\chi^{i}\chi^{j}B_{ij}(x) + \chi^{i}\chi^{j}\chi^{k} C_{ijk}(x).
\]
Recall that we have $N$ branes on $M$ with a flat $U(1)$-bundle on each. The functions under study here may depend on in which copy the string endpoints land. Thus the coefficients should be $N\times N$ matrix valued. In fact, open strings are oriented and reversing the orientation changes the sign and charge conjugates the ends which leads to coefficients in the Lie algebra of $U(N)$. 

On these constant modes it is straightforward to check that the integration functional correspond to taking trace and integrating forms over $M$, $\int\mapsto \int_{M}\tr$, that the $\star$-product becomes the usual exterior product on forms $\star\mapsto \wedge$, and from the commutation relations, the BRST-operator correspond to the exterior derivative on forms: $\mathcal{Q}\mapsto d$. A general (ghost number) argument involving formal dimensions shows that $\mathcal{A}_{0}$ is indeed a 1-form (has ghost number 1). Thus we see that the zero-mode approximation to string field theory is indeed exactly Chern-Simons theory. 

Recall now the $t$-dependence of the string action and that, as the action is $\mathcal{Q}$-exact, the theory is independent of $t$. Analyzing the $t$-dependence of the string field theory one finds that the higher oscillations can not contribute: only constant maps contributes to the path integral. Therefore the 0-mode expansion is in fact exact and we conclude that
\[
Z_{\mathrm{GW}}(T^{\ast} M,N\cdot M)\left(N,g_s=\frac{2\pi i}{k+N}\right)=Z_{\mathrm{CS}}(M)(N,k).
\]       

One can also motivate this relation more directly by observing that the open string theory localizes on holomorphic curves in $T^{\ast}M$ with boundary in $M$. There are however no non-constant such curves, so at first glance the invariant would have to vanish. The reason for it not to vanish is that the corresponding moduli space of open curves is non-compact, and for small perturbations there are solutions on very thin fat-graphs. In order to calculate the partition functions one must therefore compute the contribution from limit solutions at infinity corresponding to infinitesimally fat geodesic graphs. As indicated in our treatment of Chern-Simons perturbation theory via fat graphs, this count of degenerate Riemann surfaces agrees with the Chern-Simons partition function. This argument also says that in the general case where a Lagrangian sits inside a Calabi-Yau threefold one would expect to see instanton corrected versions of the above result with non-constant holomorphic curves on the Lagrangian giving the instanton corrections. We will see examples of this below.   

\subsection{Including knots in Chern-Simons}\label{s:knotsinCS}
So far we have discussed Chern-Simons theory and open strings for bare three-manifolds. Witten explains in \cite{Witten:1988hf} how to include knots in Chern-Simons theory: a knot is a Wilson loop observable and inserting it in the path integral gives the HOMFLY knot invariant. We will review this construction here, restricting attention to the case $M=S^{3}$.

\begin{figure}[ht]
\centering
\includegraphics[width=.6\linewidth, bb= 0 280 810 670, clip]{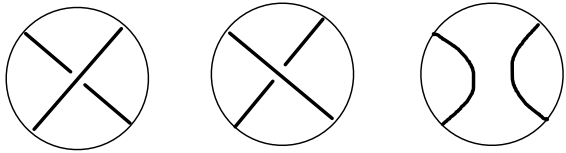}
\caption{The three knot diagrams associated to re-gluing of a small ball around a crossing in a given knot diagram.}
\label{fig:skein}
\end{figure}

Let $K\subset S^{3}$ be an oriented knot and let $A$ be a $U(N)$ connection. Write $U_{A}(K)\in U(N)$ for the holonomy of $A$ around $K$. The Wilson loop expectation value of $K$ is then
\[
W(K)=\frac{1}{Z_{\mathrm{CS}}(S^{3})}\int \mathcal{D}A \, \tr(U_{A}(K))e^{\frac{ik}{4\pi}S_{\mathrm{CS}}(A)}.
\] 
In fact for in order to make sense of this path integral we need to equip $K$ with a framing, and if none other is specified we use the homotopically unique framing for which the self linking number vanishes. Witten showed that
\[
W(K)=\frac{q^{N/2}-q^{-N/2}}{q^{1/2}-q^{-1/2}}P_{K}(q,q^{N})
\]    
where $q=\exp\left(\frac{2\pi i}{k+N}\right)$. The proof is of non-perturbative nature. It derives from a cut and paste formula for the partition function. More specifically, if a three manifold with a knot is cut into two pieces by a surface. The path integral over the two halves determine two vectors in a Hilbert space that the theory associates to the surface. Here in particular, the surface corresponds to the phase space which have finite volume. By the uncertainty relation each state occupies a phase space volume proportional to $\hbar$ and therefore the corresponding Hilbert space is finite dimensional. In the case under study a central role is played by the Hilbert space of the sphere with four marked points which is 2-dimensional. 
Consider now a knot diagram and fix a small ball around one of its crossings, see Figure \ref{fig:skein}.

 The boundary of this ball is then a sphere with four marked points. The fact that the Hilbert space is 2-dimensional then implies that there is a linear relation between the Wilson loop observables of the three links that appear from gluing the ball back in different ways. This corresponds exactly to the HOMFLY skein relation and via a small number of explicit calculations the values of the parameters involved are determined.

\subsection{Including knots in topological string theory}
The counterpart of this inclusion of knots in the open string theory was first studied by Ooguri and Vafa \cite{OV}. They associated the Lagrangian conormal $L_K\subset T^{\ast}S^{3}$:
\[
L_K=\{(q,p)\in T^{\ast}S^{3}\colon q\in K, p|_{T_qK}=0\}\approx_{\text{top}} S^{1}\times\R^{2},
\] 
to a knot $K$ and considered $M$ branes wrapping it giving an open string theory with $M$ branes wrapping $L_K$ and $N$ branes wrapping the zero section $S^{3}$. Here we restrict attention to the case $M=1$. The open string theory then has three sectors corresponding to strings with both ends on the 0-section, both ends on $L_K$, or one end on each. Arguing as above we find that the first two sectors are described by $U(N)$ Chern-Simons theory on $S^{3}$ and $U(1)$ Chern-Simons theory on $L_K$. We next study the sector of mixed strings. 

Arguing as in the string-field theory above we find that the contributions to the mixed strings come from constant strings in the intersection $S^{3}\cap L_K=S^{3}$ and that there is a single commuting $0$-mode $q$ along $K$ and a single anti-commuting $0$-mode $\chi$ with 
\[
\mathcal{A}=c(q)+A(q)\chi,
\] 
where $c(q)$ is a complex scalar. Unlike in the previous case the ambient manifold is now $T^{\ast}S^{1}$ and the formal dimension argument says that it is the scalar field that contributes so that the action for the field $c(q)$ is simply 
\[
\int_{K}\bar c dc.
\]  
However, there are gauge fields (Chan-Paton factors) $A$ on $S^{3}$ and $A'$ on $L_K$ at the ends of the string that couples to $c$ in the standard way:
\[
\int_K\tr(\bar c A c- c A'\bar c).
\]
Integrating out the $c$-field gives the following determinant 
\[
\exp\left(-\log \det\left(\frac{d}{ds}+(A_i-A'_i)\frac{dq_i}{ds}\right)\right),
\]
which can be expressed as
\[
\det(1-e^{-x}U_A(K))^{-1},
\]
where $e^{x}$ is the holonomy of the $U(1)$-connection along the generator of $H_{1}(L_K)$. To relate this to the HOMFLY invariants we expand the determinant as
\[
\det(1-e^{-x}U_{A}(K))= \sum_{S_k}\tr_{S_k} U_{A}(K) e^{-kx},
\]     
where the sum ranges over the $k$-dimensional symmetric representations of $U(N)$. Thus computing in $U(N)$ Chern-Simons theory the following weighted sum of expectation values
\begin{equation}\label{eq:waveHOMFLY}
\Psi_{K}(x)=\sum_{k} \langle \tr_{S_k} U_{A}(K)\rangle e^{-kx},
\end{equation}
which, in analogy with the discussion in Section \ref{s:knotsinCS}, is expressed through HOMFLY polynomials, and hence is polynomial in $q=e^{g_s}$ and $q^{N}$,
gives the topological string partition function for $N$ branes on $S^{3}$ and one on $L_K$:
\begin{equation}\label{eq:wavestring}
\Psi_{K}(x)= Z_{\mathrm{GW}}(T^{\ast} S^{3}, N\cdot S^{3},L_K)/Z_{\mathrm{GW}}(T^{\ast} S^{3}, N\cdot S^{3}).
\end{equation}

\begin{rmk}
When comparing the current text with \cite{AENV}. It may be useful to note that $U(N)$ Chern-Simons theory factors into $SU(N)$ and $U(1)$ Chern-Simons theories.
\end{rmk}

\section{Large $N$ transition and mirrors of the resolved conifold}\label{S:largeN}
In this section we discuss the relation between the large $N$ limit open string theory in $T^{\ast}S^{3}$ with $N$ branes on $S^{3}$ and closed string theory in the resolved conifold found by Gopakumar and Vafa, see \cite{GV}. We discuss how to move Lagrangian conormals through the transition and show how this leads to another A-model string theory. The section ends with the derivation of a local parameterization of the polynomial curve which is one of the incarnations of the main object of study in this note and which should give a mirror of the resolved conifold.

\subsection{Closed topological strings in the conifold}
Recall the free energy for the open topological string in $T^{\ast} S^{3}$ with $N$ branes on $S^{3}$:
\[
F_{\mathrm{GW}}(T^{\ast} S^{3},N\cdot S^{3})=\sum_{g,h} K_{g,h} N^{h}g_s^{2g+h-2}.
\]
For $t=Ng_s$, if we let $C_{g}(t)=\sum_{h} K_{g,h} t^{h}$ then 
\[
F_{\mathrm{GW}}(T^{\ast} S^{3},N\cdot S^{3})=\sum_{g} C_{g}(t) g_s^{2g-2},
\]   
which looks like the free energy of a closed string theory parameterized by $t$. Gopakumar and Vafa found a setup so that this is indeed the case as follows.

Consider the cotangent bundle $X=T^{\ast}S^{3}$. As a symplectic manifold its ideal contact boundary is the unit cotangent bundle $\pa_{\infty} X=ST^{\ast} S^{3}$ with the contact form given by the restriction of the action form $p\,dq$. Here we consider $T^{\ast}S^{3}$ as the algebraic variety $X_{\epsilon}$ in $\C^{4}$ with coordinates $(w_1,\dots,w_4)$ given by the equation
\[
w_1^{2}+w_2^{2}+w_3^{2}+w_4^{2}=\epsilon,\quad \epsilon\in[0,\infty),
\] 
as $\epsilon\to 0$ the $0$-section in $X_{\epsilon}$ shrinks and at $\epsilon=0$ we find a cone which at infinity looks just like $T^{\ast}S^{3}$ itself. In this way we may consider $X_{\epsilon}$ as a resolution of the cone
\[
w_1^{2}+w_2^{2}+w_3^{2}+w_4^{2}=(w_1+iw_2)(w_1-iw_2)-(iw_3+w_4)(iw_3-w_4)=0.
\]
There is another resolution $Y_{t}$: let $[\xi:\eta]$ be homogeneous coordinates on $\C P^{1}$ and consider
\begin{align*}
Y=\{&([\xi:\eta],w_1,w_2,w_3,w_4)\in\C P^{1}\times\C^{4}\colon\\
&\eta(w_1+iw_2)=\xi(iw_3-w_4),\; \eta(iw_3+w_4)=\xi(w_1-iw_2)\}.
\end{align*}
Then $Y$ is the total space of the bundle $\mathcal{O}(-1)\oplus\mathcal{O}(-1)\to\C P^{1}$. Furthermore, the projection $\pi\colon Y\to\C^{4}$ is an isomorphism outside $\pi^{-1}(0)$ and $\pi^{-1}(0)=\C P^{1}$. We consider a family $Y_{t}$ of such resolutions parameterized by the area of the sphere, concretely for $Y_t$, we take $(\xi,\eta)$ in the $3$-sphere $|\xi|^{2}+|\eta|^{2}= t$ and let $t\to 0$ then we get a 1-parameter family $Y_t$ of resolutions parameterized by $t$, which is proportional to the area of $\C P^{1}$ in the naturally induced metric, that converges to the cone at $t=0$.

\begin{figure}[ht]
\labellist
\small
\pinlabel \text{in $X_\epsilon$, $\epsilon\to 0$}  at 370 685
\pinlabel \text{in $Y_t$, $t>0$}   at 375 350
\endlabellist
\centering
\includegraphics[width=.6\linewidth]{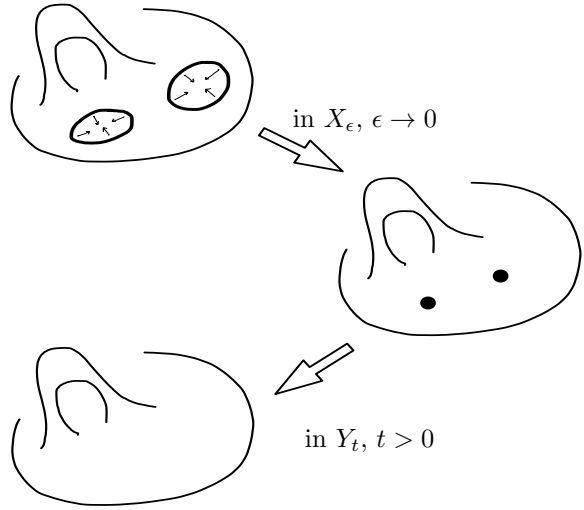}
\caption{A schematic picture of a holomorphic curve changing under conifold transition.}
\label{fig:condensation}
\end{figure}

We can now state Gopakumar-Vafa's result with more precision, it relates the open string partition function on $X$ with $N$ branes along $S^{3}$ with the closed string partition function $Y_t$ as follows:
\[
Z_{\mathrm{GW}}(X,N\cdot S^{3},g_s)=Z_{\mathrm{GW}}(Y_{t=Ng_s},g_{s}).
\]  
This conjecture has been checked via both mathematical and physical direct calculations (where the left hand side has been computed in Chern-Simons theory). There are also a number of physics proofs, see \cite{Ooguri_Vafa_worldsheet} for a proof in line with the intuition that, for the described choice of parameters as the $0$-section collapses, the boundaries of the holomorphic curves on the 0-section shrinks in an orderly manner and the curves reappear as closed curves on the other side of the transition as indicated in Figure \ref{fig:condensation}.

\subsection{Conifold transition and conormals of knots}\label{s:transitionandknots}
We next consider the Lagrangian conormal $L_K$ associated to a knot. Note that it projects to $K\subset S^{3}$. Consider a small tubular neighborhood $\approx S^{1}\times D^{2}$ of $K$. Let $d\theta$ be the standard closed $1$-form that generates $H^{1}(S^{1})$. The shift of $L_K$ along $d\theta$ is still Lagrangian and furthermore is disjoint from the zero-section. Thus we can view $L_K\subset Y_t$ since the conifold transition can be taken to affect only a small neighborhood of the vertex of the cone. In light of the intuitive explanation of the mechanism behind the equality of the partition functions it is clear that the corresponding conjecture should be as follows in this case.

Consider the open topological string theory in $X$ with $N$ branes along $S^{3}$ and one along $L_K$. The shifted $L_K$ is far from the vertex of the cone and the parts of the holomorphic curves that lie on $L_K$ should not be affected by the transition, the boundary components on the zero-section disappears in an orderly manner just as for the closed string. This leads to the conjecture
\begin{equation}\label{eq:waveconifold}
\Psi_{K}(x)=Z_{\mathrm{GW}}(Y_{t},L_K,x)/Z_{\mathrm{GW}}(Y_{t}),
\end{equation}  
where exactly as above $t=Ng_s$.

We point out that there is nothing very special about the conormal Lagrangian $L_K$ here. One would expect that the same relation holds for more general Lagrangians $L$ that can be moved off the zero-section.  

\subsection{Mirrors of the resolved conifold}\label{s:mirror}
Consider next the open string theory in $Y$ with a Lagrangian brane on $L_K$ (shifted off the zero-section). Consider the sector of strings with both endpoints on $L_K$. Here there is the contribution from short strings that lie entirely in a small neighborhood of $L_K$ which are controlled by string field theory, as above. In the present case there are non-constant holomorphic curves with boundary on $L_K$ that will give quantum corrections. The contribution from the short strings can be computed via $GL(1)$ Abelian Chern-Simons theory. The $GL(1)$ gauge group corresponds to the complex valued field $r+i\rho$, where the imaginary part is the flat connection on the $U(1)$-bundle on $L_K$ and where the real part measures the shift along $d\theta$. 

Recall that $L_{K}\approx S^{1}\times\R^{2}$. We use coordinates $(\xi,\eta,t)\in S^{1}\times S^{1}\times[0,\infty)$ on the solid torus, where $(t,\eta)$ are polar coordinates on $\R^{2}$. The connection 1-form is then the complex valued form
\[
A=A_\xi d\xi + A_\eta d\eta + A_t dt,
\]
using parallel translation in the $t$-direction we find a trivialization such that $A_t=0$ ($A_t=0$ is our gauge fixing condition). The Chern-Simons action is then, with $\dot A$ denoting the $t$-derivative of $A$,
\begin{align*}
\int_{T^{2}\times\R} A\wedge dA &= \int_{\R}dt\int_{T^{2}}d\xi\wedge d\eta\left( \dot A_{\xi} A_{\eta} -\dot A_{\eta} A_{\xi}\right)\\
& = \int_{\R} (p\dot x - x\dot p )dt=
\end{align*}
where 
\[
x=\int_{\lambda} A,\quad p=\int_{\mu} A
\]
are the periods of the connection, $\lambda=S^{1}\times\mathrm{pt}$ and $\mu=\mathrm{pt}\times S^{1}$. Note that this is the usual classical action for the form $pdx-xdp$ and hence the path integral for the Chern-Simons action is just the path integral for quantum mechanics in dimension one in terms of the periods. The usual quantization scheme then says that $x$ and $p$ should be viewed as operators with commutator
\[
[x,p]=g_s.
\] 
In the spirit of string field theory discussed above, these operators act on the partition function for strings with endpoints on $L_K$, see \eqref{eq:waveconifold}. Thus
\[
p\Psi_K(x)=g_s\frac{\partial}{\partial x}\Psi_{K}(x),
\]
and $\Psi_K(x)$ can be regarded as a wave function in quantum mechanics. In particular, the usual short wave asymptotic formula holds:
\[
\Psi_{K}(x)=\exp\left(\frac{1}{g_s}\int pdx + \dots\right),
\] 
where $\dots$ denotes terms that are $\mathcal{O}(1)$ in $g_s$. 

On the other hand $\Psi_K(x)$ is expressed \eqref{eq:waveconifold} as a Gromov-Witten partition function:
\[
\Psi_{K}(x)=\exp\left(\frac{1}{g_s}W_K(x,Q)+\dots\right),
\]
where $W_K(x)=\sum_{k,m}a_{k,m} e^{kx}Q^{m}$ is the Gromov-Witten disk potential, i.e.~the part of the free energy that counts holomorphic disk configurations of holomorphic curves with boundary on $L_K$. More complicated holomorphic curve configurations of genus $g$ and with $h$ boundary components contribute only to the coefficient of $g_{s}^{2g-2+h}$, $2g-2+h\ge 0$ of the free energy. 

From these two equations we conclude that, at the level of the disk,
\[
p=\frac{\partial W_K}{\partial x}(x,Q).
\]

Recall \eqref{eq:waveHOMFLY} the relation between the wave function $\Psi_K(x)$ and HOMFLY invariants of $K$. Using known recursive properties of the HOMFLY, see \cite{Garoufalidis}, we find that this equation, being the classical limit or characteristic variety of a holonomic D-module, is in fact (part of) an algebraic curve $V_K$ given by
\[
A_K(e^{x},e^{p},Q)=0.
\]  
Since the potential $W_{K}(x,Q)$ sums over all disk instantons on $L_K$ we can interpret the curve as the quantum corrected moduli space of the Lagrangian $L_K$. In analogy with SYZ-mirror constructions it was conjectured in \cite{AV} that, regarding $Q$ as a fixed parameter, the 3-dimensional complex hypersurface $X_{K}$ given by the equation
\[
A_{K}(e^{x},e^{p},Q)=uv
\] 
is a B-model mirror to the resolved conifold $Y$. Mirror symmetry exchanges A-branes, Lagrangian submanifolds, in $Y$ with B-branes, holomorphic submanifolds, in $X_K$. In particular, associated to every Lagrangian brane $L$ in $Y$ there is a B-brane in $X_{K}$ so that the quantum corrected moduli space of $L$ in $Y$ is the same as the classical moduli space of the mirror B-brane in $X_K$. Here the B-brane mirror to $L_K$ at $(e^{x_0},e^{p_0})$ is the line $\{u=0, v\text{ arbitrary}\}$ over the point $(e^{x_0},e^{p_0})\in V_K$.

\section{Knot contact homology}\label{S:kch}
Knot contact homology was introduced and studied by Ng, see \cite{Ngframed, Ngtransverse, Ngsurvey}. It associates a differential graded algebra (DGA) to any knot or link in $S^{3}$. Here we will describe this theory from the point of view of its contact geometric origins: it is the Legendrian contact homology algebra of the unit conormal lift of a knot. Legendrian contact homology is a Floer type theory and in particular uses moduli spaces of holomorphic disks. In the case under study the DGA admits a combinatorial description: it can be calculated from a braid presentation of a knot. Here we will not say much about concrete calculations but rather focus on the geometric features of the theory that allows us to relate it to the topological string theory described above. Legendrian contact homology is a part of Symplectic field theory introduced by Eliashberg, Givental, and Hofer \cite{EGH} that was developed in a setting relevant to knot contact homology by Ekholm, Etnyre, Ng, and Sullivan in a series of papers \cite{EES,E,EENS, EENStransverse}. 

\subsection{Contact geometry and knots}
Let $K\subset M$ be a knot in a three-manifold $M$. As mentioned above we view the cotangent bundle $T^{\ast}M$ as an exact symplectic manifold with Liouville form $\beta=pdq$ and symplectic form $\omega=d\beta$. The ideal contact boundary of $T^{\ast}M$ is the 5-dimensional spherical cotangent bundle $ST^{\ast}M$ with contact form $\alpha=\beta|_{ST^{\ast}M}$, the contact condition is $\alpha\wedge d\alpha^{2}\ne 0$. The contact form $\alpha$ determines a Reeb vector field $R_{\alpha}$, by $d\alpha(R_{\alpha},\cdot)=0$ and $\alpha(R_{\alpha})=1$. Closed trajectories of $R_{\alpha}$ will be called Reeb orbits. Representing the spherical cotangent bundle as the unit cotangent bundle using some Riemannian metric the Reeb flow is exactly the geodesic flow on the unit cotangent bundle and Reeb orbits correspond to closed geodesics. 

The ideal boundary of the conormal Lagrangian $L_K\subset T^{\ast}M$ is a Legendrian submanifold $\Lambda_{K}\subset ST^{\ast}M$, which topologically is a torus that can be naturally identified with the boundary of small tubular neighborhoods of the knot. In particular, the homology $H_{1}(\Lambda_{K})$ comes equipped with a canonical basis $\lambda,\mu$, where $\lambda$ is the longitude and $\mu$ the meridian of the knot $K$. Trajectories of $R_{\alpha}$ starting and ending on $\Lambda_{K}$ will be called Reeb chords. In the unit conormal bundle Reeb chords correspond to bi-normal chords, i.e.~geodesic curves with endpoints on the link and intersecting it at right angles at these endpoints.

\subsection{Legendrian contact homology}
Contact homology is a part of Symplectic field theory and is defined for general contact manifolds with Legendrian submanifolds in them. Here we will however describe it only in the case relevant to knot contact homology. As with topological strings contact homology has a closed and an open sector. The closed sector is a commutative DGA $\mathcal{B}=\mathcal{B}(ST^{\ast}M)$ generated by its Reeb orbits and the open sector is a non-commutative DGA $\mathcal{A}=\mathcal{A}(\Lambda_K)$ generated by Reeb chords and with coefficients in $\mathcal{B}$ that is called Legendrian contact homology. Here we will be mainly interested in the case $M=S^{3}$ in which case the closed and open sectors decouple and the chord algebra $\mathcal{A}$ can be studied separately without considering orbits. For simplicity, we will describe the construction of the algebra $\mathcal{A}$ only in this case, and in order not to have to consider orbits we take the link under study to lie in a small ball $\approx\R^{3}$ and one can show, compare \cite{EkholmNg}, that the theory reduces to the corresponding theory in $\R^{3}$, in other words, we can restrict attention to the short Reeb orbits lying over this ball since long chords leaving the ball and orbits have no effect on the algebra. 

Let $\mathcal{A}=\mathcal{A}(\Lambda_{K})$ denote the algebra over the group ring $\C[H_{2}(ST^{\ast}S^{3},\Lambda_K)]$ generated by the Reeb chords of $\Lambda_{K}$. We associate a grading $|c|$ to a chord $c$ via a certain Maslov index. In the case under study this index is simply the usual Morse index of the corresponding bi-normal chord which takes the values $0,1,2$. We pick generators $t,x,p$ of $H_{2}(ST^{\ast}S^{3},\Lambda_K)$, where $t$ is the class of the fiber $2$-sphere and where $x$ and $p$ map to the longitude
and the meridian in $H_1(\Lambda_{K})$. We then write the coefficients of $\mathcal{A}$ as 
\[
\C[H_{2}(ST^{\ast}S^{3},\Lambda_K)]=\C[e^{\pm x},e^{\pm p}, Q^{\pm 1}],\quad Q=e^{t}.
\]
In order to keep track of homology classes of certain maps we fix some auxiliary data that we describe next. For each Reeb chord $c$ we fix a capping disk $\gamma_{c}\colon D\to ST^{\ast}S^{3}$, where $D$ is the unit disk in the complex plane. Here, $\gamma_{c}$ restricted to the boundary arc between $-1$ and $1$ parameterizes the Reeb chord $c$ and it maps the complementary arc into $\Lambda_K$. 

\begin{figure}[ht]
\labellist
\small
\pinlabel $a$  at 225 780
\pinlabel $b_1$ at 30 60
\pinlabel $b_2$ at 152 60
\pinlabel $b_m$ at 467 60
\endlabellist
\centering
\includegraphics[width=.5\linewidth]{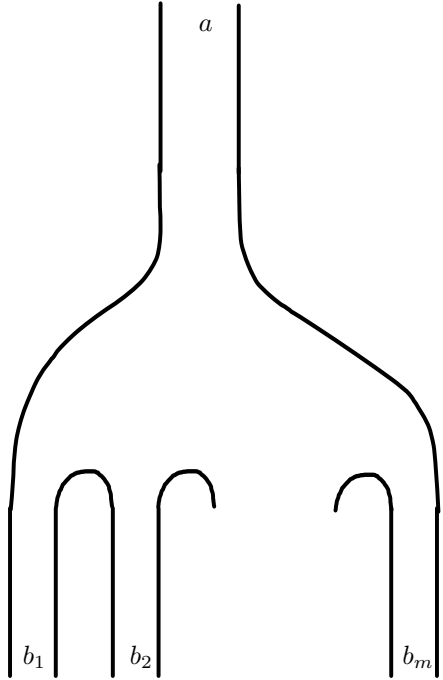}
\caption{A holomorphic disk in $\R\times ST^{\ast} S^{3}$ with boundary on $\R\times\Lambda_{K}$ and Reeb chord asymptotics.}
\label{fig:disk1}
\end{figure}

The differential on $\mathcal{A}$ is defined via a count of holomorphic disks. We start with a general discussion of the moduli spaces of disks that we will use. Consider the symplectization $\R\times ST^{\ast}S^{3}$ with the symplectic form $d(e^{s}\alpha)$, where $s$ is a coordinate in the $\R$-factor. (Note that this is symplectomorphic to the complement of the zero-section in the cotangent bundle, $T^{\ast}S^{3}-S^{3}$.) Pick an almost complex structure $J$ on $\R\times ST^{\ast}S^{3}$ that decomposes as a complex structure in the contact planes compatible with $d\alpha$ and such that $J\partial_{s}=R_{\alpha}$. In this almost complex structure, if $c$ is a Reeb chord then $\R\times c$ is a holomorphic strip with boundary on $\R\times\Lambda_K$. Let $D_{m}$ be a a disk with one distinguished boundary puncture, called positive, $m$ other boundary punctures, called negative. Consider maps $u\colon (D_{m+1},\partial D_{m+1})\to(\R\times  ST^{\ast}S^{3},\R\times\Lambda_{K})$ which are asymptotic to a holomorphic strip $\R\times a$ at $+\infty$ at the positive puncture and asymptotic to the holomorphic strip $\R\times b_j$ at $-\infty$ at the $j^{\rm th}$ negative puncture. Write $\mathbf{b}=b_1\dots b_m$ and define $\mathcal{M}_{e^{nx},e^{kp},Q^{l}}(a;\mathbf{b})$ to be the moduli space of holomorphic such disk maps $u$, i.e.~$du+J\,du\,i=0$, such that when $u$ is capped off  by the capping disks at the Reeb chords $a,b_1,\dots,b_m$ then it represents the homology class $nx+kp+lt$, see Figure \ref{fig:disk1}.

Then its formal dimension is
\[
\dim(\mathcal{M}_{e^{nx}e^{kp}Q^{l}}(a;\mathbf{b}))=|a|-|\mathbf{b}|.
\]   
We point out that moduli spaces of this form are empty unless the action of the Reeb chord at the positive puncture is at least as large as the sum of the actions of the Reeb chords at the negative punctures. This is a consequence of Stokes theorem as follows: for a holomorphic disk $u\in\mathcal{M}_{e^{nx}e^{kp}Q^{l}}(a;\mathbf{b})$, $u^{\ast}(d\alpha)$ is positive since $J$ is compatible with $d\alpha$  and thus 
\[
\int_{a}\alpha-\sum_j\int_{b_j}\alpha=\int_{D_m}u^{\ast}(d\alpha)\ge 0.
\]

\begin{figure}[ht]
\labellist
\small
\pinlabel $\dim=d_1$  at 170 610
\pinlabel $\dim=d_2$ at 85 215
\pinlabel $\dim=d_3$ at 250 215
\endlabellist
\centering
\includegraphics[width=.7\linewidth]{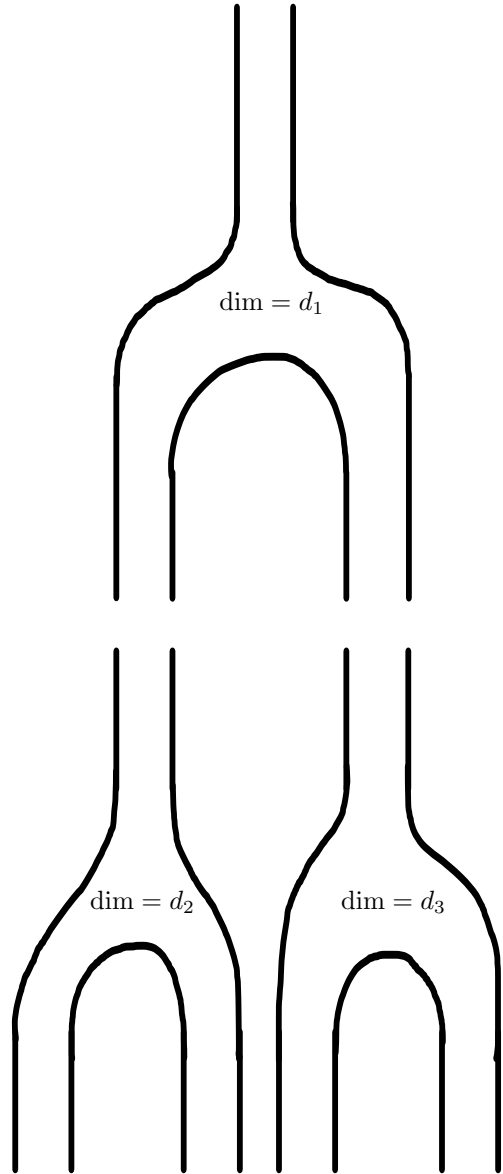}
\caption{A two level holomorphic disk in the boundary of a moduli space of dimension $d_1+d_2+d_3$.}
\label{fig:twolevel}
\end{figure}

By our choice of almost complex structure, $\R$ acts on solutions of the holomorphic disk equation. It can be shown that for generic $J$ the moduli spaces considered above are transversely cut out. (Transversality is particularly simple in this case because the condition of having only one positive punctures rules out multiple covers.) The moduli spaces $\mathcal{M}_{e^{nx}e^{kp}Q^{l}}(a,\mathbf{b})$ admit a nice compactification (after dividing by the $\R$-action), analogous to breaking of flow lines in Morse theory, where the boundary of a moduli space consists of several level disks joined at Reeb chords, see \cite{BEHWZ}. Furthermore, from the dimension formula above the total dimension of a building is the sum of the dimensions of its levels, see Figure \ref{fig:twolevel}.

Fixing a spin structure on $\Lambda_{K}$ we get an induced orientation on the index bundle associated to the linearization of the holomorphic disk equation of the configuration space of maps $u$. This then induces an orientation on the transversely cut out moduli space. These orientations of the index bundles can be chosen coherently, or in other words to be compatible with gluing. It then follows in particular that if $\mathcal{M}_{e^{nx}e^{kp}Q^{l}}(a,\mathbf{b})$ is a moduli space of dimension $2$ then after dividing by the $\R$-action we get a reduced oriented moduli space of dimension $1$ and the boundary of this moduli space consists of two level buildings where each level is a curve of dimension $1$ and where the product of the orientation signs of the corresponding reduced moduli spaces at the ends of the 1-manifold agrees with the sign induced from the orientation of the 1-manifold, see Figure \ref{fig:twolevel}.

The differential $\partial\colon\mathcal{A}(\Lambda_{K})\to\mathcal{A}(\Lambda_{K})$ satisfies Leibniz rule and is defined through the following holomorphic disk count on generators:
\[
\pa a = \sum_{\begin{smallmatrix}|a|-|\mathbf{b}|=1\\ n,k,l\in\Z\end{smallmatrix}}
|\mathcal{M}_{e^{nx}e^{kp}Q^{l}}(a,\mathbf{b})|\,e^{nx}e^{kp}Q^{l}\,\mathbf{b},
\]    
where $|\mathcal{M}|$ denotes the algebraic number of points in the reduced moduli space $\mathcal{M}/\R$, which is an oriented 0-manifold. By the above mentioned compactness result the reduced moduli space is compact and thus the sum above is finite. To see that $\partial^{2}=0$ we note that by the above description of 2-dimensional moduli spaces, two level disks which contributes to $\partial^{2}a$ are in oriented one to one correspondence with the boundary of an oriented 1-manifold and hence cancel out.

\subsection{Computations}

\begin{ex}
We compute the knot contact homology of the unknot using a heuristic argument. For full details see \cite{EENS}. Consider the unknot $O$ as a round planar circle. Note that there is an $S^{1}$-family of bi-normal chords of index $1$. The Reeb chords of $\Lambda_{O}$ thus come in an $S^{1}$ Bott-family of index $1$. Since the action of all chords are equal the only possible disks with only one positive puncture has no negative punctures. From a variational perspective, the holomorphic disk equation can be thought of the gradient flow equation for the action functional $\gamma\mapsto \int_{\gamma}\alpha$ for curves $\gamma$ with endpoints on $\Lambda_{K}$. A chord of $O$ can flow in two directions and then shrinks until its length is $0$, see Fig \ref{fig:circleinplane}. In order to control what happens for lengths near zero we look at the path traced out by the endpoints of the flowing chord in $\Lambda_{O}$, see Figure \ref{fig:circleinplane}.
\begin{figure}[ht]
\centering
\includegraphics[width=.6\linewidth]{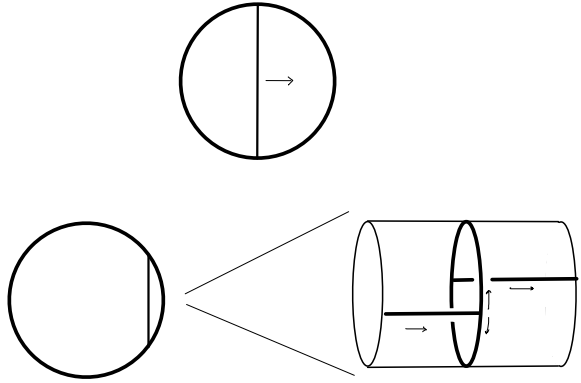}
\caption{Top: start of chord flow. Bottom: the path traced out by the endpoints of the chord flow near chord length zero on the conormal torus is closed up by the boundary values of the two holomorphic disks.}
\label{fig:circleinplane}
\end{figure}

There are two ways to close the curve up and using a concrete perturbation of $\Lambda_{O}$ one can show that there is exactly one holomorphic disk in each class. In conclusion there are four disks with positive puncture at any chord in the Bott family. Morisfication of the Bott situation gives two chords $e$ and $c$ of gradings $|e|=2$ and $|c|=1$. The differential of $e$ corresponds to the Morse differential in the circle:
\[
\partial e = c-c=0,
\]
and the differential of $c$ gets contributions from the four disks just discussed. It is clear that the boundaries of these disks lie in distinct homology classes in $H_{1}(\Lambda_O)$. In order to find the component of the homology class along the fiber sphere we count intersections with an $S^{3}$-section that  agrees with covector dual to the direction of propagation of one of the chord flows. We find that the differential is
\[
\partial c = 1-e^{x}-e^{p}+Qe^{x}e^{p}.
\] 
Although there are four signs in this formula three of them can be changed by various choices but then the last is fixed and is computed by a determinant bundle calculation. 
\end{ex}

With this example established we give a short description of how knot contact homology is calculated for general links. We start by representing the link $K$ as a braid around the unknot $O$. As we push the braid toward the unknot, $\Lambda_{K}$ approaches $\Lambda_{O}$ and in particular if the braid is sufficiently close to the unknot, $\Lambda_{K}$ is a graphical Legendrian in a standard contact neighborhood of $\Lambda_{O}$ which is contactomorphic to the $1$-jet space of the tours $J^{1}(\Lambda_O)$. We then show that in the limit $K\to O$ rigid holomorphic disks on $\Lambda_{K}$ are in natural one to one correspondence with rigid configurations of holomorphic disks on $\Lambda_{O}$, that are controlled by the above description of the disks in the Bott-degerate case, and Morse flow trees of $\Lambda_{K}\subset J^{1}(\Lambda_{O})$, that can be described from the braid presentation. 
 
\subsection{Exact cobordisms, DGA-morphisms, and invariance} 
Contact homology has functorial properties similar to TQFTs. Here we describe them only in the setup where we need them and refer to \cite{EGH,E} for more general considerations. Consider first two knots $K$ and $K'$ with corresponding Legendrian tori $\Lambda_{K}$ and $\Lambda_{K'}$. Assume that there exists an exact Lagrangian cobordism $L\subset \R\times ST^{\ast}S^{3}$ such that $L$ agrees with $[M,\infty)\times\Lambda_{K}$ for $s>M$ and with $(-\infty,-M]\times\Lambda_{K'}$ for $s<-M$. Here the exactness condition means that the form $e^{s}\alpha|_{L}$, which is closed since $L$ is Lagrangian, is exact, $e^{s}\alpha=dz$ for some $z\colon L\to\R$. Note that there is an inclusion map $H_{2}(ST^{\ast}S^{3},\Lambda_{K})\to H_{2}(\R\times ST^{\ast}S^{3},L)$ and similar for the negative end. Using these maps we change the ground ring of both algebras $\mathcal{A}(\Lambda_{K})$ and $\mathcal{A}(\Lambda_{K'})$ to $\C[H_{2}(\R\times ST^{\ast}S^{3},L)]$ 

In analogy with the above we write $\mathcal{M}^{L}_{A}(a,\mathbf{b})$ for the moduli space of holomorphic disks $u\colon (D_{m+1},\partial D_{m+1})\to(\R\times ST^{\ast}S^{3},L)$ with positive puncture asymptotic to the Reeb chord strip $\R\times a$ for a Reeb chord of $\Lambda_{K}$ in the positive end and to Reeb chord strips $\R\times b_j$ in the negative end, where $\mathbf{b}=b_1\dots b_m$ is a word of Reeb chords of $\Lambda_{K'}$, which when closed up by capping disks represents the homology class $A\in H_{2}(\R\times ST^{\ast}S^{3},L)$.

Define the DGA map $\Phi_{L}\colon\mathcal{A}(\Lambda_{K})\to\mathcal{A}(\Lambda_{K'})$ to act on generators as
\[
\Phi_{L}(a)=\sum_{\begin{smallmatrix}|a|-|\mathbf{b}|=0\\ A\in H_2(\R\times ST^{\ast}S^{3},L)\end{smallmatrix}} |\mathcal{M}^{L}_{A}(a,\mathbf{b})|e^{A}\mathbf{b},
\]
where this time $|\mathcal{M}^{L}|$ denotes the algebraic number of points in the oriented 0-manifold $\mathcal{M}^{L}$. In analogy with our argument above showing that $\partial^{2}=0$ we find that $\Phi_{L}$ is a chain map, $\Phi\partial-\partial\Phi=0$. To see this we simply note again that the two level disks that contributes to the terms in the chain map equation (i.e.~with one level of dimension 0 in the cobordism and one level of dimension one in the positive or negative end) are in oriented one to one correspondence with the ends of a compact  oriented 1-manifold. Similarly, if $L'$ is a cobordism with positive end $\Lambda_{K'}$ and negative end $\Lambda_{K''}$ then we can glue the two cobordisms by the level $s=R$ of $L'$ to the level $s=R$ of $L$ giving a cobordism $L''_R$. Here we consider a 1-parameter family of cobordisms for $R\to\infty$ and also in this case the moduli space admits a compactification with several level disks at the boundary. In particular for $R$ sufficiently large rigid disks in the cobordism are in oriented one to one correspondence with two level disks where both levels are rigid. This shows that for $R$ sufficiently large
\[
\Phi_{L''_{R}}=\Phi_{L'}\circ\Phi_{L}.
\]     

A similar argument, that however uses a more advanced perturbation scheme, see \cite{Ekholm_rsft}, shows that the homotopy class of map $\Phi_L$ is invariant under deformations (e.g.~of $L$ or the almost complex structure). This implies that $\mathcal{A}(\Lambda_K)$ is invariant under Legendrian isotopies of $\Lambda_{K}$ as follows. A Legendrian isotopy connecting  $\Lambda_{K}$ to $\Lambda'_{K}$ gives cobordisms $L_0$ and $L_1$ connecting $\Lambda_{K}$ to $\Lambda_{K'}$ and vice versa. Composing these cobordism we get a cobordism that can be deformed to the trivial cobordism which clearly induces the identity map on $\mathcal{A}(\Lambda_K)$ or $\mathcal{A}(\Lambda_{K}')$. We thus conclude that
\[
\Phi_{L_0}\circ\Phi_{L_1}\simeq 1,\quad \Phi_{L_1}\circ\Phi_{L_0}\simeq 1,
\]     
which proves that $\Phi_{L_0}$ and $\Phi_{L_1}$ are homotopy inverses.

Below we will also make use of another type of exact Lagrangian cobordisms. Consider an exact Lagrangian $L\subset T^{\ast}S^{3}$ with ideal boundary $\Lambda_{K}$. Then $L$ has empty negative end and the maps above are thus simply maps into the algebra of the empty Lagrangian which is simply the ground ring $\C[H_{2}(T^{\ast}S^{3},L)]$ with the trivial differential. This map counts holomorphic disks with one positive puncture and boundary on $L$:
\begin{equation}\label{e:exactchmap}
\Phi_{L}(a)=\sum_{\begin{smallmatrix}|a|=0\\ A\in H_{2}(T^{\ast}S^{3},L)\end{smallmatrix}}|\mathcal{M}_{A}^{L}(a)|\,e^{A}.
\end{equation}
Note that, since the differential in the target is trivial in this case, the chain map equation simply reads $\Phi_L\circ\partial=0$. This type of map is an example of an augmentation and will play a central role in all that follows. 

\begin{rmk}
We point out that the exactness condition on $L$ is crucial for the chain map properties discussed above. More precisely, the exactness condition shows that any closed holomorphic disk $D$ with boundary on $L$ would have area $\int_{D}d(e^{s}\alpha)=\int_{\partial D}e^{s}\alpha=\int_{\partial D}dz=0$ and hence would have to be constant. This rules out boundary bubbling and shows that the moduli spaces have compactifications with several level disks joined at Reeb chords only.
\end{rmk}
\section{Augmentations and the augmentation variety}\label{S:aug}
We will next study augmentations of $\mathcal{A}(\Lambda_{K})$ more systematically. Recall that the coefficient ring of $\mathcal{A}(\Lambda_{K})$ was the Laurent polynomial ring $\C[e^{\pm x},e^{\pm p},Q^{\pm 1}]$. Assigning non-zero complex values to our formal variables, i.e.~taking $(e^{x},e^{p},Q)\in(\C^{\ast})^{3}$, we think of $\mathcal{A}(\Lambda_{K})$ as a family over $(\C^{\ast})^{3}$ of algebras over $\C$. We then define an augmentation of $\mathcal{A}(\Lambda_{K})$ to be a DGA-map (i.e.~a chain map which is the canonical map on coefficients): 
\[
\epsilon\colon \mathcal{A}(\Lambda_{K})\to\C,\quad \epsilon\circ\partial-\partial\circ\epsilon=\epsilon\circ\partial=0,
\]    
where we think of $\C$ as a unital DGA in degree zero with the trivial differential.

We define the \emph{augmentation variety} $V_K$ of $\mathcal{A}(\Lambda_K)$ to be the closure of the highest dimensional part of the locus of points $(e^{x},e^{p},Q)\in(\C^{\ast})^{3}$ where $\mathcal{A}(\Lambda_{K})$ admits augmentation. In all examples this is a codimension one subvariety. Its defining polynomial $A_{K}(e^{x},e^{p},Q)$ is the augmentation polynomial.

\begin{ex}
We compute the augmentation polynomial of the unknot. Here the condition that $\epsilon$  is an augmentation is simply
\[
\epsilon(\pa c)=\epsilon(1-e^{x}-e^{p}+Qe^{x}e^{p})=0
\] 
and thus the augmentation polynomial is simply:
\[
A_{O}=1-e^{x}-e^{p}+Qe^{x}e^{p}.
\]
\end{ex} 

\subsection{Augmentations and exact Lagrangian fillings}
Note that for $Q=1$, the augmentation polynomial for the unknot reduces to
\[
A_{O}(e^{x},e^{p},Q=1)=(1-e^{x})(1-e^{p}).
\]
This is no coincidence, in fact the lines $\{Q=1,e^{p}=1\}$ and $\{Q=1,e^{x}=1\}$ lies in the augmentation variety for any knot $K$. For the first line, we note that the conormal $L_K$ gives an exact Lagrangian filling of $\Lambda_{K}$ in $T^{\ast}M$, and that in the map $H_{2}(ST^{\ast} S^{3},\Lambda_K)\to H_{2}(T^{\ast}S^{3},L_{K})$ takes $p$ and $t$ to zero (recall $Q=e^{t}$) and takes $x$ to $x$, viewed as the generator of $H_{1}(L_K)$. For the second line we note that $L_K\cap S^{3}=K$ is a clean intersection of Lagrangians. Performing Lagrangian surgery along $K$ we construct an exact Lagrangian filling $M_K$ of $\Lambda_K$ which has the topology of the knot complement $S^{3}-K$. In this case the map $H_{2}(ST^{\ast} S^{3},\Lambda_K)\to H_{2}(T^{\ast}S^{3},M_{K})$ takes $x$ and $t$ to zero and takes $p$ to $p$, viewed as the generator of $H_{1}(M_K)$. 

\subsection{Augmentations and non-exact Lagrangian fillings}\label{s:nonexact}
Recall from Section \ref{s:transitionandknots} that we shifted $L_K$ off of the zero-section in $T^{\ast} S^{3}$ and that after conifold transition we can view it as a Lagrangian $L_K\subset Y$. Note however that the Lagrangian submanifold $L_K$ is no longer exact. In particular, as the physics arguments of Section \ref{s:transitionandknots} shows there are closed holomorphic disks (i.e.~without punctures) with boundary on $L_{K}$. As we shall see, these holomorphic disks interferes with the definition of an augmentation that we used in the exact case. In fact, the resolution of that problem, using obstruction and bounding cochains in the spirit of \cite{FO3} will allow us to relate the augmentation polynomial and the polynomial defining the B-model mirror with B-branes corresponding to the knot Lagrangians.     
 
We first note that the shifted copy of $L_K$ shifts $\Lambda_{K}$ along the closed form $d\theta$ at infinity. In \cite{AENV} we argued that although $\Lambda_{K}$ was shifted there is a natural correspondence between punctured holomorphic disks at infinity with boundary on the original and the shifted  $\Lambda_{K}$. Furthermore, for closed disks below any given area limit, there exists a compact subset of $Y$ that contains all of them. We next explain the role of boundary bubbling in defining chain maps. 

Consider a moduli space $\mathcal{M}^{L_{K}}(a)$ of holomorphic disks with a positive puncture at the Reeb chord $a$ and boundary on $L_{K}$. Also, in this non-exact case there is a version of Gromov compactness, the limit is a several level holomorphic building, with holomorphic disks in the symplectization at infinity, but now the level in $Y$ may itself be a broken disk consisting of one disk with a positive puncture with several holomorphic disks attached along its boundary. Here the boundary bubbling is analogous to what happens to the real curve in $\C^{2}$ with boundary on $\R^{2}$ given by the equation $xy=\epsilon$ as $\epsilon\to 0$.

\begin{figure}[ht]
\centering
\includegraphics[width=.8\linewidth, bb= 0 225 790 610, clip]{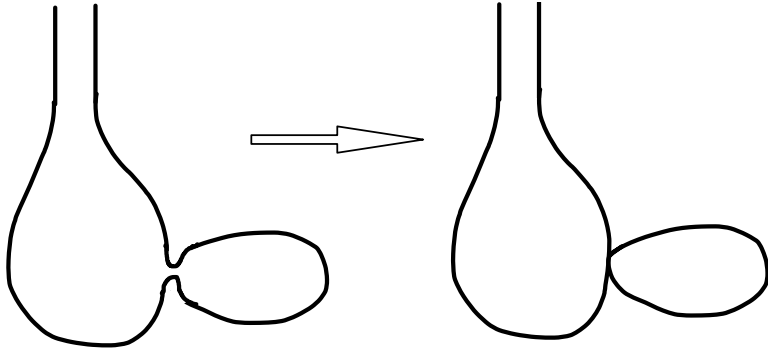}
\caption{Boundary bubbling introduces new boundary components.}
\label{fig:boundarysplit}
\end{figure}

In particular, in the case that the dimension of $\mathcal{M}^{L_K}(a)$ is one we find that the boundary may contain broken disks with one component in $\mathcal{M}^{L_K}(a)$ and one rigid disk on $L_K$. This means that the chain map equation for the map $\Phi_L$ in \eqref{e:exactchmap} does not hold, see Figure \ref{fig:boundarysplit}. We will next describe how to modify the definition of that map to resolve this problem.

Assume that we have fixed a perturbation scheme so that all these disk (below a cut off area) are transversely cut out. In order to resolve the problem with boundary bubbling we use a construction similar to the bounding cochain techniques used in \cite{FO3}. Pick a reference curve $\xi$ of class $x$ in $\Lambda_{K}$ and then pick for each rigid disk  $u\colon (D,\pa D)\to(Y,L_K)$ an obstruction killing 2-chain $\sigma_u$ in $L_K$ that connects $u(\partial D)$ with $k\xi$, where $kx$ is the homology class of $u(\partial D)$. Let $\sigma$ denote the collection of all obstruction killing 2-chains and define the moduli space $\mathcal{M}_{e^{kx}Q^{l}}^{L_K}(a;\sigma)$ to be the space of holomorphic disks with positive puncture at $a$ and with all possible insertions of $\sigma$ along its boundary. More precisely, an element in this moduli space is a several level object with finitely many levels as follows. The top level is a holomorphic map $v\colon D\to (Y,L_K)$, with one boundary puncture and finitely many boundary marked points, where the disk is asymptotic to the Reeb chord strip of $a$ at the positive puncture and where the evaluation map at each marked point takes value in some $\sigma_{u}$. The disks $u$ corresponding to the obstruction killing chains constitute the second level and are also equipped with additional marked points that map to obstruction chains of other disks that then constitutes the third level and so on, attaching the obstruction chains of the disks in the $(r+1)^{\rm th}$ level at marked points on the boundary of the disks in the $r^{\rm th}$ level.  

We define the homology class of such a configuration as (recall $Q=e^{t}$)
\[
kx+lt=(k_0x+l_0t)+(k_1x+l_1t)+\dots+(k_mx+l_mt),
\]    
where $k_0x + l_0t$ is the class of the disk $v$ and where $(k_jx+l_jt)$ are the homology classes of the rigid holomorphic disks in the building and we let $\mathcal{M}_{e^{kx}Q^{l}}^{L_K}(a;\sigma)$ denote the moduli space of configurations of total homology class $kx+lt$ (but with no restriction on the specific number of insertions).

We define the map $\epsilon_{L_K}\colon\mathcal{A}(\Lambda_K)\to\C[e^{\pm x},Q^{\pm 1}]$ as follows
\[
\epsilon_{L_{K}}(a)=\sum_{\begin{smallmatrix}|a|=0\\ k,l\in\Z\end{smallmatrix}}
|\mathcal{M}_{e^{kx}Q^{l}}^{L_K}(a;\sigma)|\,e^{kx}Q^{l},
\]
where $|\mathcal{M}|$ denotes the algebraic number of configurations. 

We next define the GW-potential as the generating function of the holomorphic disk configurations that are inserted in a disk with insertions. More precisely, let $u_1$ and $u_2$ be holomorphic disks with boundary on $L_K$ and let $\sigma_{u_1}$ and $\sigma_{u_2}$ be the corresponding obstruction chains. Then $\sigma_{u_1}$ and $\sigma_{u_2}$ both agree with a multiple of a standard curve $\xi$ at infinity. Shifting $\xi$ off itself we make the two chains disjoint outside a compact subset and after small perturbation the chains are transverse elsewhere. Define the linking number of $u_1$ and $u_2$ as the intersection number between $u_{1}(\partial D)$ and $\sigma_{u_2}$ or equivalently as the intersection number of $u_2(\partial D)$ and $\sigma_{u_1}$. The objects that we count are now finite weighted trees with holomorphic disks at the vertices. Here the weight at a vertex is the weight of the rigid disk as a component of the moduli space as weighted manifold, the weight of an edge is the linking number of the two disks at its endpoints and the total weight of the tree is the product of weights at all its edges and vertices. We write $\mathcal{M}(L_K,\sigma)$ for the moduli space of such tree configurations and we define the homology class of a tree configuration as the sum of the homology classes of the disks in the tree. 

The GW-potential then has the form:
\[
W_K(x,Q)=\sum_{k,r\in\Z} C_{k,r}\,e^{kx}Q^{r},
\]
where $C_{k,r}$ is the sum of the weights of the tree configurations in $\mathcal{M}(L_K;\sigma)$ that represent the homology class $kx+rt$.

Our main result can now be formulated as follows: if $p=\frac{\partial W_{K}}{\partial x}(x,Q)$ then the map $\epsilon_{L_{K}}\colon\mathcal{A}(\Lambda_{K})\to\C[e^{\pm x},Q^{\pm 1}]$ is a chain map. It then follows that
\[
p=\frac{\partial W_K}{\partial x}(x,Q)
\]  
gives a local parameterization of the augmentation variety $V_K$. In particular this variety agrees locally with the B-mirror curve defined in \ref{s:mirror} and the corresponding polynomials must agree in the irreducible case.

\begin{figure}[ht]
\centering
\includegraphics[width=.8\linewidth, bb= 0 270 790 610, clip]{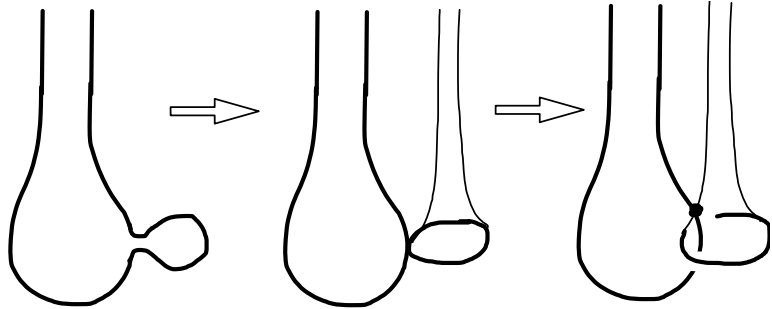}
\caption{Obstruction killing 2-chains turn boundary splittings into interior points in the moduli space.}
\label{fig:obstrchain}
\end{figure}

\begin{figure}[ht]
\labellist
\small
\pinlabel \text{in $\R\times ST^{\ast} S^{3}$} at 500 575
\pinlabel \text{in $Y_t$} at 500 210
\endlabellist
\centering
\includegraphics[width=.4\linewidth]{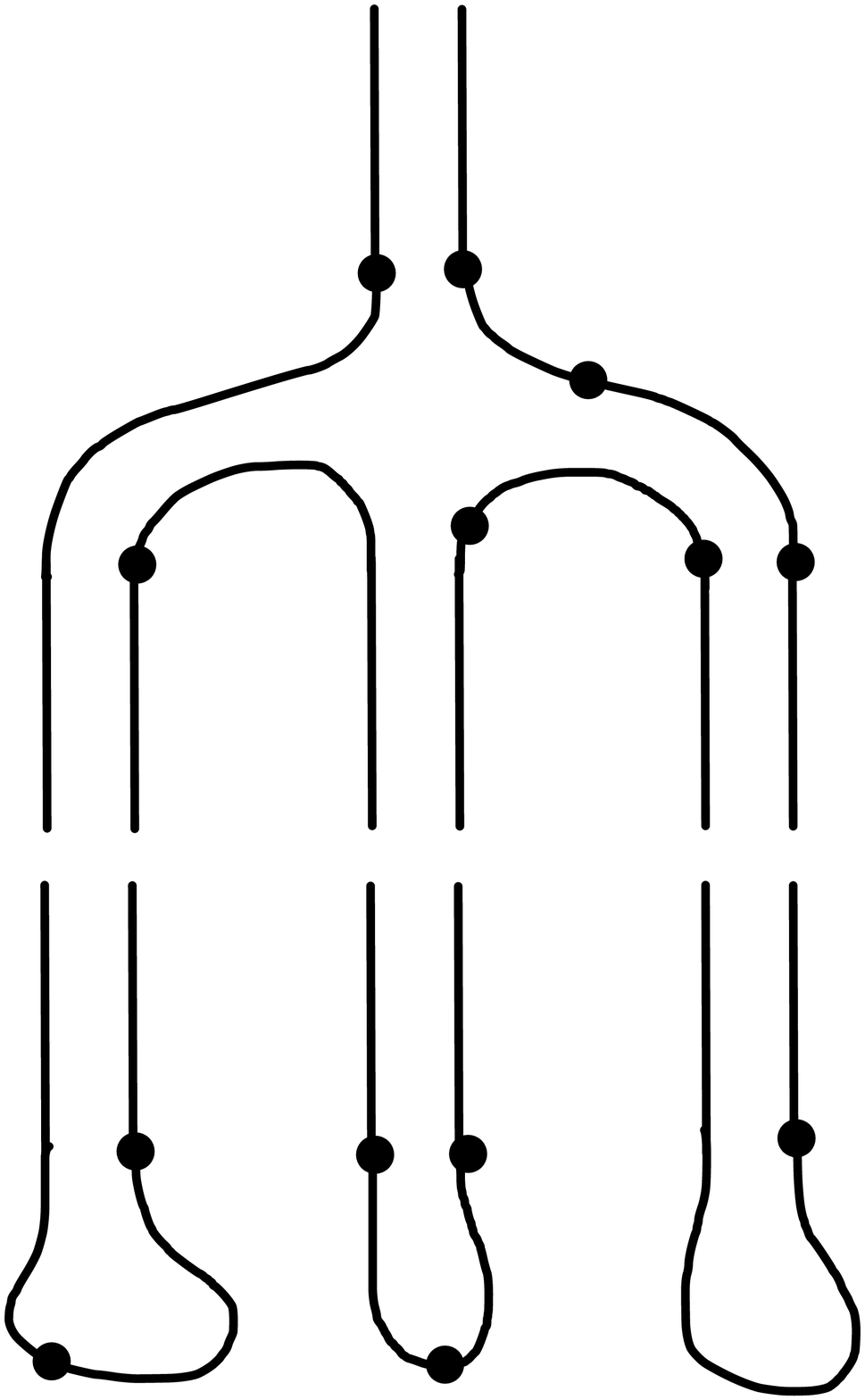}
\caption{Boundary configuration in the space of holomorphic disks with insertions. The dots indicates inserted obstruction killing 2-chains.}
\label{fig:insertions}
\end{figure}

To show this result we study the boundary of the 1-dimensional moduli space of disks with insertions. First observe that the obstruction killing 2-chain turns boundary splitting into interior points of the moduli space, as indicated in Figure \ref{fig:obstrchain}.
Consequently, the boundary of a 1-dimensional space of disks with insertions consists of two level disks with insertions. The lower level is a rigid disk with insertions and correspond exactly to the configurations contributing to the maps $\epsilon_{L_{K}}$ The upper level is a 1-dimensional family of disks with insertions in the symplectization end, see Figure \ref{fig:insertions}.
By construction, the obstruction killing chain looks like a product of the $\R$-factor and a multiple of the reference curve $\xi$ in $\Lambda_{K}$ in homology class $x$. This allows us to derive an algebraic expression for the count of disks with insertions in the symplectization end. 

Let $v\colon D_{m+1}\to\R\times ST^{\ast}S^{3}$ be a holomorphic disk in $\mathcal{M}_{e^{nx}e^{kp}Q^{l}}(a;\mathbf{b})$. Consider the disjoint union of all the obstruction killing 2-chains in $\sigma$ with the corresponding configuration rigid holomorphic disks attached at lower end and, by construction, looking like $\R\times \xi^{r}$ in the symplectization end, where $\xi^{r}$ is a multiple of the reference curve $\xi$ in $\Lambda_{K}$. Then the Gromov-Witten potential $W_{K}(e^{\pm x},Q)$ is the generating function counting the elements in this disjoint union. To see this note that the linking number between two rigid disks counts the number of possible insertions.

Consider next the disjoint union of obstruction killing 2-chains decorated with a marked point where they can be inserted in the disk $u$. Such a decorated 2-chain consists of a 2-chain $\sigma_u$ asymptotic to $\R\times \xi^{r}$ together with a point in $\xi^{r}$ at which $\sigma_{u}$ is attached to $v(\pa D_{m+1})$. If $\sigma_{u}$ is asymptotic to $\R\times\xi^{r}$ such points come in $r$-tuples since $\xi^{r}$ has multiplicity $r$. Furthermore, any chain with a marked point correspond to an intersection point between $\xi$ and $v(\pa D_{m+1})$ and as such come equipped with an intersection sign, counting them with signs we find that there are $k$ such points (recall that the homology class of $v(\pa D_{m+1})$ is $nx+kp$). In conclusion the obstruction killing 2-chains $\sigma_{u'}$ of a rigid disk configuration $u'$ in homology class $rx+ht$ gives rise to $kr$ chains decorated with a marked point where they can be can be attached. It then follows that
\[
k\frac{\partial W_{K}}{\partial x}= \sum krM_{r,h}e^{rx}Q^{h} 
\]       
is the generating function counting such 2-chains with marked points that can be inserted. Consider next
\begin{align*}
e^{k\frac{\partial W_{K}}{\partial x}} &=1+ \left(\sum krM_{r,h}e^{rx}Q^{h}\right)\\
& + \frac{1}{2}\left(\sum krM_{r,h}e^{rx}Q^{h}\right)^{2}+\dots+\frac{1}{m!}\left(\sum krM_{r,h}e^{rx}Q^{h}\right)^{m}+\dots
\end{align*}
and note that the first term $1$ counts no insertion and that the general $m^{\rm th}$ term counts $m$ insertions. It follows that
\[
e^{nx}e^{k\frac{\partial W_K}{\partial x}}Q^{l}
\] 
counts the disks with insertions arising from one disk in $\mathcal{M}_{e^{nx}e^{kp}Q^{l}}(a,\mathbf{b})$. This then implies that for $p=\frac{\partial W_{K}}{\partial x}$, $\epsilon_{L_K}$ has the chain map property and thus is an augmentation which is the desired result.  
   
The above argument is not special to the conormal. A similar argument shows that any Lagrangian which looks like a shift of $\Lambda_{K}$ along some closed form at infinity and which has first homology generated by the first homology of $\Lambda_{K}$ induces augmentations via its Gromov-Witten disk potential in a similar way. In fact, as explained in \cite{AENV} one can also use immersed Lagrangian fillings to construct augmentations provided they are unobstructed, which means that the count of rigid disks with one Lagrangian corner with insertions vanishes. It is a very interesting problem to find constructions of such Lagrangians, in particular since they could help proving that the augmentation polynomial gives the B-model mirror in general.

\subsection{A brief discussion of the case of links}\label{s:links}
It is straightforward to generalize the discussion of knot contact homology above to the case of many  component links $K=K_1\cup\dots\cup K_n$. Here the conormal $\Lambda_{K}$ is a many component Legendrian with torus components. The coefficient ring for the knot contact homology algebra is accordingly $\C[e^{\pm x_j},e^{\pm p_j},Q^{\pm 1}]$, where $x_j$ and $p_j$ are the longitude and meridian of the $j^{\rm th}$ knot component $K_{j}$. In this case the augmentation variety (for fixed $Q$) is a subvariety in $(\C^{\ast})^{2n}$. This variety has several components. The simplest of them is the one corresponding to augmentations of the individual torus components of $\Lambda_{K}$. These give augmentations of the link that map any mixed chord, i.e.~Reeb chords with endpoints on distinct torus components of the Legendrian link to zero. The corresponding component of the augmentation variety is then just the product of the augmentation curves. 

One can also consider augmentations that are non-trivial on mixed chords. Conjectural constructions give unobstructed connected Lagrangian fillings $M_K\subset Y$ of $\Lambda_{K}$ and, much like in the knot case, there is a disk super-potential $W_{K}(p_1,\dots,p_n)$ associated to them that gives a local parameterization of the augmentation variety as
\[
x_j=\frac{\partial W_{K}}{\partial p_j},
\]      
which is Lagrangian in $(\C^{\ast})^{2n}$.

In \cite{AENV} we conjecture that $V_{K}$ is always Lagrangian and that it is the characteristic variety of a holonomic  $\mathcal{D}$-module with generator the HOMFLY wave function $\Psi_{K}(x_1\dots,x_n)$ of the link. Furthermore, we conjecture that this $\mathcal{D}$-module arises from A-model open topological strings in $(\C^{\ast})^{2n}$ with a space filling coisotropic brane and a Lagrangian brane on $V_K$. 
Here the algebra of open string states with both endpoints on the coisotropic brane leads to viewing the coordinates on $(\C^{\ast})^{2n}$ as operators satisfying the commutation relations of the Weyl algebra and $\Psi_{K}(x_1,\dots,x_n)$ is the wave function of a fermion in $V_K$ that corresponds to the unique ground state of strings with one endpoint on the coisotropic brane and one on $V_K$. We refer to \cite{AENV} for more information on this conjecture.    
 
\bibliographystyle{utcaps}
\bibliography{myrefs}

\end{document}